\documentclass{article}
\usepackage{graphicx}
\usepackage{amsthm}

\newcommand{\nc}{\newcommand}
\nc{\rnc}{\renewcommand}

\nc{\al} {\alpha}
\nc{\eps}{\epsilon}
\nc{\G}{\Gamma}
\rnc{\l}{\lambda}
\nc{\sig}{\sigma}
\rnc{\t}{\tau}
\rnc{\th}{\theta}
\nc{\w}{\omega}

\nc{\C}{\ensuremath{{\bf C}}}
\rnc{\d}{{\rm d}}
\nc{\Dsp}{\displaystyle}
\nc{\e}{\ensuremath{{\rm e}}}
\nc{\fr}{\frac}
\rnc{\i}{\ensuremath{{\rm i}}}
\nc{\limK}{\lim\limits_{K\to\infty}}
\nc{\lb}{\left[}
\nc{\lp}{\left(}
\nc{\lsb}{\left\{}
\nc{\nn}{\nonumber}
\nc{\noin}{\noindent}
\nc{\N}{\ensuremath{{\bf N}}}
\rnc{\O}[1]{\ensuremath{O\left( #1 \right)}}
\nc{\rb}{\right]}
\nc{\rp}{\right)}
\nc{\rsb}{\right\}}
\nc{\R}{\ensuremath{{\bf R}}}
\nc{\Z}{\ensuremath{{\bf Z}}}

\nc{\ba}{\begin{array}}
\nc{\ea}{\end{array}}
\nc{\be}{\begin{equation}}
\nc{\ee}{\end{equation}}
\nc{\bea}{\begin{eqnarray}}
\nc{\eea}{\end{eqnarray}}
\nc{\beas}{\begin{eqnarray*}}
\nc{\eeas}{\end{eqnarray*}}

\nc{\intRp}{\int_0^{\infty}}

\newcounter{theorem}
\newtheorem{thm}[theorem]{Theorem}
\newtheorem{cor}[theorem]{Corollary}

\newtheorem{lemma}[theorem]{Lemma}

\nc{\alK}{\lfloor \al K \rfloor}
\nc{\Arg}{{\rm Arg}}
\nc{\betaK}{\lfloor K/\al \rfloor}
\nc{\cosPiS}{\cos\lp\fr{\pi}{2}s\rp}
\nc{\dl}{\d\l}
\nc{\fOne}{\ensuremath{f_1(\t,\l)}}
\nc{\fTwo}{\ensuremath{f_2(\t,\l)}}
\rnc{\L}{\ensuremath{\Lambda}}
\nc{\LK}{\ensuremath{\L^K}}
\nc{\LKab}[2]{\ensuremath{\L^K_{(#1,#2)}}}
\nc{\LO}{\ensuremath{\Lambda\setminus\{0\}}}
\rnc{\P}{\ensuremath{\wp}}
\nc{\ReSTwo}{\ensuremath{\Re(s)\!>\!2}}
\nc{\sumE}{\lim_{N\to\infty}\sum_{n=-N}^N\lim_{M\to\infty}\sum_{m=-M}^M}
\nc{\sumHD}{\lim_{K\to\infty}\sum_{|m|,|n|\le K}}
\nc{\sumLat} {\sum_{\w\in\LO}}
\nc{\sumTay} {\sum_{n=3}^{\infty}}
\nc{\sumNotE}{\lim_{M\to\infty}\sum_{m=-M}^M\lim_{N\to\infty}\sum_{n=-N}^N}
\nc{\sumQuadK}[2]{\sum_{#1=1}^K\sum_{#2=-#1}^{#1}\!\!\eps_{#1#2}}
\nc{\sumLKab}[2]{\sum_{\w\in\LKab{#1}{#2}}\!\!\!}
\nc{\tE}{\widetilde{E}}
\nc{\tEs}{\ensuremath{\tE_s}}
\nc{\tETwoA}{\ensuremath{\tE_2^{(\al)}}}
\nc{\wmn}[2]{\w_{#1,#2}}
\nc{\z}{\ensuremath{\zeta}}

\begin{document}
\title{Integral Representations for Elliptic Functions} 
\author{Andrew Dienstfrey\footnote{National Institute of Standards and
Technology, Boulder, CO 80305-3328.}\ and Jingfang 
Huang\footnote{University of North Carolina, Chapel Hill, NC 27599-6011.}}

\maketitle 

\begin{abstract}
We derive new integral representations for objects arising in the
classical theory of elliptic functions: the Eisenstein series $E_s$,
and Weierstrass' $\wp$ and $\zeta$ functions.  The derivations proceed
from the Laplace-Mellin transformation for multipoles, and an
elementary lemma on the summation of 2D geometric series.  In
addition, we present new results concerning the analytic continuation
of the Eisenstein series as an entire function in $s$, and the value
of the conditionally convergent series, denoted by $\widetilde{E}_2$
below, as a function of summation over increasingly large rectangles
with arbitrary fixed aspect ratio\footnote{Contribution of
U.S. Government, not subject to copyright.}.
\end{abstract}

\section{Introduction}
In this paper we revisit the classical theory of elliptic functions as
developed by Eisenstein and Weierstrass.  Both of these researchers
represented the meromorphic functions appearing in their theories as
summations over a given lattice of elementary pole functions of a
prescribed order.  Our fundamental observation is that pole functions
may be represented by exponentially-damped, oscillatory integrals.
These representations depend on the complex half-planes in which the
singularities lie, and are natural variants of the classical Mellin,
or Laplace-Mellin, formulas which are valid for isolated poles lying
in the right half plane, where $\Re(\w)>0$.  In more recent times such
integral representations have resurfaced in the development of fast
multipole methods for the solution of Poisson problems; in this
context they often are referred to as ``plane-wave'' representations
(\cite{hr98} and, more recently, \cite{cgr99}).  A key feature of
these integral representations is that the pole centers appear in the
exponents of the integrands.  As a consequence the lattice summations
are transformed into geometric series which may be summed explicitly
underneath the integral.  The result is a new class of integral
representations for the Eisenstein series and other meromorphic
functions of Weierstrass' theory.

A brief summary of the paper follows.  In the first section we
review the definitions of the Eisenstein series $E_n$ and the
Weierstrass functions \P\ and \z.  We will analyze a generalization
of Eisenstein's series which we denote by $\tEs$, the differences
being: first, we consider $s=\sig+\i t\in\C$, and second, we define
$\tEs$ as a limit over lattice squares of increasing size, a
significant point when $\Re(s)\le2$ and the sums are not absolutely
convergent.  In addition, in this preliminary section we provide
elementary derivations of the requisite plane-wave formulas for
general pole functions of the form $f(\w)=\w^{-s}$, and a summation
identity for a two-dimensional geometric series.

In the next section we employ the plane-wave representations and the
summation identity to derive an integral representation for $\tEs$ for
the case $\ReSTwo$.  Integral representations for Eisenstein's $E_n$
naturally follow for $s=n\ge 3$.  Subsequently, we derive an
alternative representation for $\tEs$ as a contour integral from which
we deduce that the sums $\tEs$, defined unambiguously for \ReSTwo,
admit an analytic continuation as an {\bf entire} function to the
whole of the complex plane. As a corollary, we prove the existence of
a finite limit for $\tE_2$.  We discuss $\tE_2$ and its relation to
Eisenstein's, $E_2$.  As the summation processes defining these two
conditionally convergent series are distinct, one expects different
limiting values.  We derive a closed form correction term for this
difference.  (Note that both Eisenstein's convention and ours give
$E_1=\tE_1=0$.) In the following section we derive analogous integral
formulas for Weierstrass' \P\ and \z\ functions.  We conclude the
paper with a brief discussion of these integral representations in
relation to previous research in the theories of lattice sums, and
elliptic functions.

In this last regard we take a moment to mention here that our formulas
for $\tEs$ are the natural lattice analogues to the well-known
representation for Riemann's zeta function (which we denote with the
subscript $\z_R$ so as to distinguish it from Weierstrass' function
of the same name)
\be
\label{RiemannZ}
 \z_R(s) 
   = \sum_{n=1}^\infty \fr{1}{n^s}
   = \fr{1}{\Gamma(s)}\intRp \l^{s-1} \fr{1}{1-\e^{-\l}} \e^{-\l} \dl,
   \ \Re(s)>1.
\ee
For example, in the case of a square lattice we derive the following
integral expression for the classical Eisenstein series
\bea
   E_k(\i) 
    &=& \sum_{n=-\infty}^\infty \sum_{m=-\infty}^\infty
        \fr{1}{(m+n\i)^k} \nn \\
\label{EkSquare}
    &=& \fr{8}{(k-1)!}\intRp \l^{k-1}  
        \fr{\cos^2(\l/2)}{1-2\e^{-\l}\cos(\l)+\e^{-2\l}} \e^{-\l} \dl,
\eea
where $k$ is a positive integer divisible by four ($E_k(i)=0$
otherwise).  The similarity between (\ref{RiemannZ}) and
(\ref{EkSquare}) is clear.  For more general lattices, we replace $\i$
by $\t$, $k\in\N$ by $s\in\C$, and the single trigonometric ratio in
(\ref{EkSquare}) by a sum of analogous ratios denoted by \fOne\ and
\fTwo\ defined in (\ref{fOne:def}) and (\ref{fTwo:def}).  The general
expression is given in theorem~\ref{Eisenstein:thm}.

We note that a subset of the results presented below appeared
previously in a slightly different form \cite{huang99}.

\section{Preliminaries}
In this section we review the definitions of the Eisenstein series and
the Weierstrass \P\ and \z\ functions.  Furthermore, we derive
elementary lemmas concerning plane-wave representations and a
geometric series identity, both of which we will use repeatedly in the
subsequent sections.

\subsection{The Eisenstein series and Elliptic functions}
We are given a general lattice $\Lambda\subset \C$ which we describe
by generators $\mu, \nu$: 
\[
   \Lambda=\{m\cdot\mu + n\cdot\nu| m,n\in\Z\} 
\]
where $\mu,\nu$ are complex numbers such that the lattice ratio,
$\t=\nu/\mu$, is not real.  The fundamental parallelogram is the
set, $\L_0=\{\al+\beta\t| |\al|\le 1/2, |\beta|\le 1/2\}$.  Eisenstein
began his investigations of doubly periodic, meromorphic functions
with periods $\mu$ and $\nu$ through the study his eponymous series
(see, for example, \cite{weil76} and
\cite{schoeneberg74})
\be
\label{EisensteinEn}
 E_k = \sumE \frac{1}{(m\cdot \mu + n\cdot \nu)^{k}},\ k\ge 1;
\ee
the elimination of the term $m=n=0$ is implicit here and below. For
the magnitude of each term in the summand we have the bounds
\be
\label{MagBound}
    \fr{k(\mu,\nu)}{\sqrt{m^2+n^2}} 
\le \fr{1}{|m\cdot\mu + n\cdot\nu|} 
\le \fr{K(\mu,\nu)}{\sqrt{m^2+n^2}}.
\ee
Applying (\ref{MagBound}) and elementary estimates one verifies that
the series (\ref{EisensteinEn}) are absolutely convergent for $n\ge
3$, and thus are well-defined functions of the lattice
$E_n=E_n(\Lambda)$.  The same estimates indicate that the series $E_n$
are absolutely divergent for $n=1$ or $2$, hence the limiting
operation specified in (\ref{EisensteinEn}) plays a non-trivial role
in the definition of these sums.  Eisenstein proved that the limiting
procedure (\ref{EisensteinEn}) yields finite values of $E_n$ even in
these cases.  We observe that, in particular, $E_1=0$.  This is in
keeping with the fact that, using a symmetry argument in the
absolutely convergent case, one may prove that $E_n=0$ for all odd
$n\ge3$.  Therefore, from the point of view of convergence, the only
``interesting'' sum is $E_2$.

Eisenstein was cognizant that the value of $E_2$ depends on the choice
of limiting procedures; and he derived many identities which connect
his summation process for $E_2$ to others (\cite{weil76}). We choose
yet a different summation convention and define $\tEs$ as the limit of
partial sums over ``lattice-squares'' of increasing size.  We
generalize further in considering complex exponents.  Specifically, we
define $\tEs$ by
\be
\label{HDEn}
 \tEs = \sumHD \frac{1}{(m\cdot \mu + n\cdot \nu)^{s}},
\ee
which we consider, initially, for $s=1, s=2$ and $\Re(s)>2$.  For
non-integer $s$ we consider the branch of the function $\z^{s}$ with a
cut along the ``negative diagonal'' of the lattice, $\{z=-t(\mu+\nu),\
t>0\}$.  For all $z$ in the closure of this cut plane we have
\be
\label{BranchW}
\ba{l}
   \th\le \arg(z) \le \th+2\pi, \mbox{ where }\\
   \th=\Arg(-\mu-\nu).
\ea
\ee
We further enforce the convention that points of the lattice lying
along the diagonal are considered symmetrically,
\[
 \fr{1}{(-m\mu-m\nu)^s} = \fr{1}{2(m|\mu+\nu|)^s}
    \lp\fr{1}{\e^{\i\th s}}
     + \fr{1}{\e^{\i(\th+2\pi)s}}\rp
\]
We will return to this point later.  Finally, we note that the
large $K$ limiting convention defined in (\ref{HDEn}) is relevant only
in the cases $s=1$ or $s=2$.

We will derive integral representations for $\tEs$.  Naturally,
restricting $s$ to the positive integers, our formula yields an
integral representation for the classical Eisenstein series
$\tE_n=E_n$, for $n\ge3$.  As for the conditionally convergent series,
it is straightforward to verify that $\tE_1=0$ directly from
(\ref{HDEn}).  However, for $s=2$, it is not {\it a priori} obvious
that (\ref{HDEn}) will be finite for arbitrary $\Lambda(\mu,\nu)$.
Furthermore, if finite, the relationship between its value and
Eisenstein's $E_2$ is not clear.  The finite existence of the limit
(\ref{HDEn}) will follow as a consequence of the integral
representations for $\tEs, \Re(s)>2$.  In addition, we derive a
formula which connects our limiting value to Eisenstein's.  Even more,
we prove that $\tEs$ admits an analytic continuation to $s\in\C$ as an
entire function.

Some fifteen years on the heels of Eisenstein, in 1862 Weierstrass
commenced his study of doubly-periodic functions.  For Weierstrass,
the fundamental object was his \P\ function which he defined as
\be
\label{WeierstrassP}
  \P(x,\L) = \fr{1}{x^2} + \!\sumLat\lp\fr{1}{(x-\w)^2}-\fr{1}{\w^2}\rp.
\ee
Standard estimates demonstrate that \P\ is an analytic function of
$x$, and well-defined function of $\L$.  From (\ref{WeierstrassP}) we
see that the \P\ function has a double pole at the origin and every
lattice translate.  Weierstrass defined his \z\ function as an indefinite
integral of \P\, and developed the following summation representation
\be
\label{WeierstrassZ}
  \z(x,\L) 
    = -\int^x \P(s,\L) \d s 
    =  \fr{1}{x} + \!\sumLat\lp\fr{1}{(x-\w)}+\fr{1}{\w}+\fr{x}{\w^2}\rp.
\ee
Again, the correction terms inside the summation balance the
asymptotics of the translated poles so that the sum converges
absolutely for $x$ in a compact set containing no lattice points.  As
\z\ is meromorphic and has only a single pole within a fundamental
parallelogram, it can not be doubly periodic (see, for example,
\cite{ahlfors79}).  Nevertheless, Weierstrass derived many properties
satisfied by this function and its relatives.  The absolutely
convergent sums (\ref{WeierstrassP}) and (\ref{WeierstrassZ}) will
serve as the starting points for the derivation of the integral
representations for \P\ and \z\ below.

We conclude this section with a brief comment on freedom in the choice
of generators for $\L$.  Note that $\tEs, \P, \z$ satisfy simple
rescalings with respect to $\mu$
\be
\label{Rescalings}
\lsb \ba{rcl}
  \tEs(\mu,\nu) &=&\Dsp \fr{1}{\mu^s}\tEs(1,\t), \\
  \P(x|\mu,\nu)&=&\Dsp 
     \fr{1}{\mu^2}\P\lp \left.\fr{x}{\mu}\right| 1,\t\rp,\\
  \z(x|\mu,\nu)&=&\Dsp 
     \fr{1}{\mu}  \z\lp \left.\fr{x}{\mu}\right| 1,\t\rp,
\ea \right.
\ee
where $\t=\nu/\mu$.  

Furthermore, one observes that two pairs of generators $(\mu,\nu)$ and
$(\mu',\nu')$ will give rise to the same lattice if and only if they
satisfy a linear system of the form
\[
  \lp\ba{c} \mu' \\ \nu' \ea\rp 
     = \lp\ba{cc} a&b\\c&d \ea\rp \lp\ba{c} \mu \\ \nu \ea\rp
\]
with $a,b,c,d\in\Z, ad-bc=\pm1$. The set of all two-by-two matrices
with integral entries and determinant plus or minus one is known as
the unimodular group in two variables. In this case the lattice ratios are
transformed by the fractional linear transformation 
\[
    \t' = \fr{a+b\t}{c+d\t},
\]
known as a unimodular substitution. One may prove that up to rescaling
and unimodular substitution, any lattice ratio may be represented by a
unique $\t$ chosen from the following fundamental region
\cite{ahlfors79}
\be
\label{FundamentalT}
\lsb\ba{l}
  \Dsp -\frac{1}{2} < \Re(\t) \le \frac{1}{2},  \\
        \Im(\t) > 0 ,\\
           |\t| \ge  1, \\
           \mbox{if } |\t|= 1,\mbox{ then } \Re(\t)\ge 0. \\
   \ea \right.
\ee

In summary, without loss of generality, we restrict our analysis to
the ``inhomogeneous'' functions, which are obtained from (\ref{HDEn}),
(\ref{WeierstrassP}) and (\ref{WeierstrassZ}) by fixing $\mu=1,
\nu=\t$, and consider $\L=\L(\t)$ with $\t$ satisfying 
(\ref{FundamentalT}).  For convenience we omit the variables $\mu,\nu$
below and write, for example, $\P=\P(z,\t)$.

\subsection{Plane-wave representations and a 2D Geometric Series}

\begin{figure}[htbp]
\label{lattice:fig} 
  \begin{center} 
  \includegraphics[height=5in,width=5in]{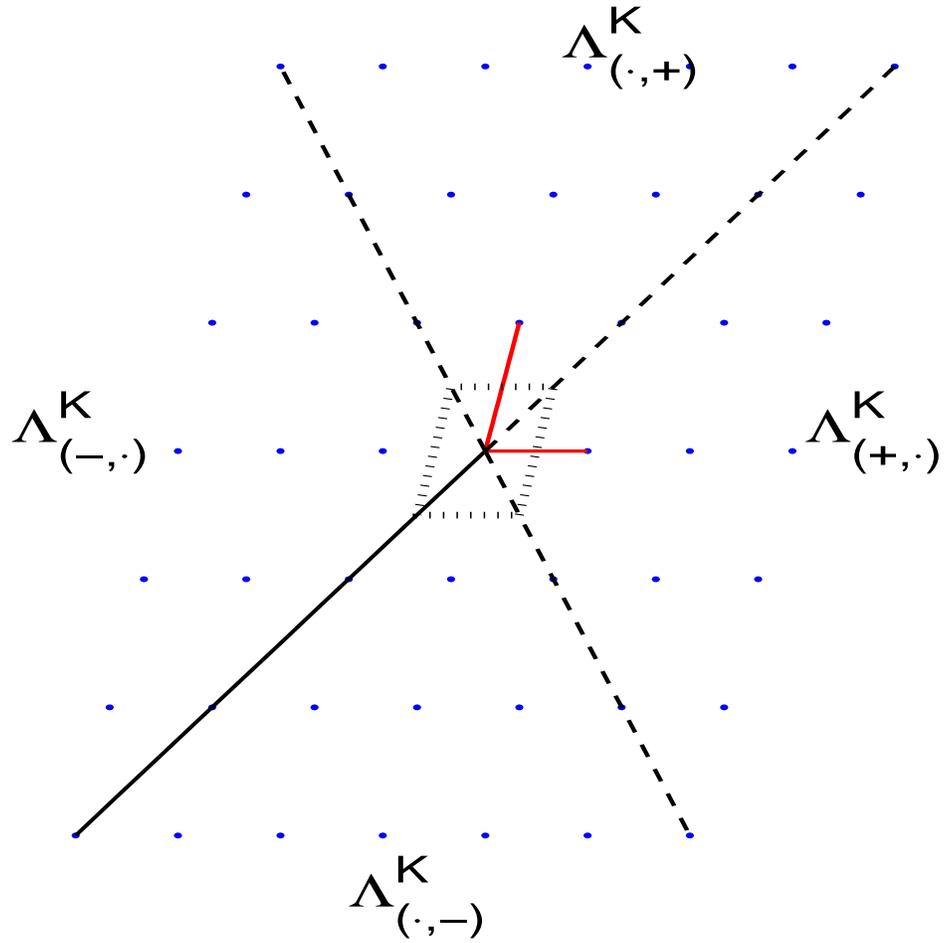}
  \caption{Partition of $\LK$ into subregions.  The generators
    $(1,\t)$ are shown in red.  The central dotted region is the
    boundary of the fundamental domain.  The dashed lines show the
    divisions into $\LKab{\pm}{\pm}$.  The solid black line is the
    branch cut.}
  \end{center}
\end{figure}

To facilitate our derivations we define the truncated lattice
$\LK=\{\wmn{m}{n}=m+n\t|\ |m|,|n|\!<\!K\}\setminus\{0\}$. We further
group lattice points into four overlapping ``quadrants''
\[
  \LK =  \LKab{+}{\bullet} \cup \LKab{\bullet}{+} 
    \cup \LKab{-}{\bullet} \cup \LKab{\bullet}{-} 
\]
defined by
\be
\addtolength{\arraycolsep}{-.1em}
\label{QuadrantDefn}
\ba{rclrcccrrcl}
 \LKab{+}{\bullet} &=& \{\wmn{m}{n}|& 1 &\le& m &\le & K,\ &|n|&\le& m\}\\
 \LKab{\bullet}{+} &=& \{\wmn{m}{n}|& 1 &\le& n &\le & K,\ &|m|&\le& n\}\\
 \LKab{-}{\bullet} &=& \{\wmn{m}{n}|&-K &\le& m &\le &-1,\ &|n|&\le& -m\}\\
 \LKab{\bullet}{-} &=& \{\wmn{m}{n}|&-K &\le& n &\le &-1,\ &|m|&\le& -n\}.
\ea
\ee
We recall from the discussion following (\ref{HDEn}) that for
non-integer $s$ the shared boundary between $\LKab{-}{\bullet}$ and
$\LKab{\bullet}{-}$ is identical to the branch cut (see
figure~\ref{lattice:fig}).

We have the following elementary lemma
\begin{lemma} 
\label{PlaneWave:lem}
Assume a complex lattice $\L(\t)$ and the quadrants defined as in
(\ref{QuadrantDefn}). An isolated singularity of complex order $s,
\ \Re(s)\!>\!0$ with branch cut defined as in (\ref{BranchW}) may be  
represented by the following plane-wave integrals, each of which is
valid in the appropriate quadrant determined by the location of the
point $\w$.
\be
\label{PlaneWave}
\fr{1}{\w^s} = \left\{ \ba{rl}
\Dsp           \fr{1}{\G(s)}\intRp \l^{s-1} \e^{- \l\w} \dl, 
  &\w\in\LKab{+}{\bullet} \vspace{1em}\\
\Dsp \fr{\e^{-\i\pi s}}{\G(s)}\intRp \l^{s-1} \e^{  \l\w} \dl, 
  &\w\in\LKab{-}{\bullet} \vspace{1em}\\
\Dsp \fr{\e^{-\i\pi s/2}}{\G(s)}\intRp \l^{s-1} \e^{\i\l\w} \dl, 
  &\w\in\LKab{\bullet}{+}  \vspace{1em}\\
\Dsp \fr{\e^{ \i\pi s/2}}{\G(s)}\intRp \l^{s-1} \e^{-\i\l\w} \dl, 
  &\w\in\LKab{\bullet}{-}.
    \ea\right.
\ee
\end{lemma}
\bigskip

\noin {\bf Proof:}\\
We have the representation of the $\G$ function
\[
  \G(s) = \intRp \l^{s-1}\e^{-\l}\dl,\ \Re(s)>0.
\]
Assume $\w\in \i\R^+$, in particular, $\w=\i t,\ t>0$. As $\t$ satisfies
(\ref{FundamentalT}) we observe that $-\pi<\th <-\pi/2$ hence
$\arg(\w^s)=\pi s/2$.  Keeping this in mind, we rescale the
integration variable by $t$, factor the $-1$ in the exponential, and
multiply and divide by $\exp(\i\pi s/2)$ to obtain
\beas
   \G(s) &=& t^{s}\intRp \l^{s-1}\e^{-\l t}\dl \\
         &=& \e^{-\i \pi s/2}\w^{s}\intRp \l^{s-1}\e^{\i\l \w}\dl.
\eeas
Dividing both sides by $\G(s)\w^s$ gives the desired result for
$\w=\i t$.  In a similar manner we prove the formula for $\w$ lying on
any of the principal coordinate rays emanating from the origin, $\w\in
\pm \R^+, \pm \i\R^+$.  The full expressions (\ref{PlaneWave}) then
follow by analytic continuation into the appropriate
quadrants. $\dagger$
\bigskip

\noin We note that for integer $s$, the integral expressions may be
continued further and are valid in the appropriate half-planes
$\pm\Re(\w)>0$ and $\pm\Im(\w)>0$.

Next we turn to our summation convention (\ref{HDEn}). From
lemma~\ref{PlaneWave:lem}, it is apparent that no single plane-wave
expansion formula will be valid for all terms in the summands
(\ref{HDEn}), (\ref{WeierstrassP}), and (\ref{WeierstrassZ}); terms
must be grouped with respect to quadrant.  As with the our convention
of splitting contributions from points $\w_{-m,-m}$ lying on the cut
in (\ref{HDEn}) equally between branches, we wish to treat each
quadrant as symmetrically as possible.  We define the symbol
$\eps_{mn}$ for $m,n\in\Z$ by
\be
\label{EpsilonSymbol:def}
\eps_{mn} = 
 \lsb \ba{rr}
      \Dsp \fr{1}{2}, & m   = \pm n \vspace{1em} \\
      \Dsp 1,         & \mbox{ otherwise}
 \ea \right.
\ee
By convention, we sum over the terms in the $\LKab{+}{\bullet}$-quadrant 
as
\bea
\label{Grouping}
\sumLKab{+}{\bullet} f(\w)
     &=& \sumQuadK{m}{n} f(\wmn{m}{n}) \\
\label{GroupingExplicit}
     &=& \sum_{m=1}^{K}
         \lp \fr{1}{2}f(\wmn{m}{-m})
            +\sum_{n=-m+1}^{m-1} f(\wmn{m}{n})
            +\fr{1}{2}f(\wmn{m}{m}) \rp
\eea
We recognize the slight abuse of notation in (\ref{Grouping}) in that
the left hand side is not interpreted as a standard sum over the set
$\LKab{+}{\bullet}$, but rather a modified sum in which the
``diagonal'' terms are added with a factor of 1/2; this is made
explicit in (\ref{GroupingExplicit}).  We choose this abuse of
notation so as to not burden our summation symbols with a cascade of
modifiers, and anticipate that it will not cause confusion below.

The sum over the lattice square \LK\ is the sum of the quadrant sums
as in (\ref{Grouping}); hence the reason for the factor of $1/2$ ---
to avoid double-counting of the contributions from the diagonal terms
--- is clear.  We note that for $s\in\N,\ s>2$ the numerical value of
the sums is independent of any manner of grouping terms.  Even so, the
form of the integrands in our integral representations reflect this
choice.  We have found that the convention (\ref{Grouping}) yields the
most symmetric expressions in appearance (a different grouping for
integer $s$ was employed in \cite{huang99}).

Next, we recognize that in sums of the form (\ref{Grouping}) we may
substitute the appropriate plane-wave expansion (\ref{PlaneWave}) to
represent the poles contained in $f$.  This substitution will
transform the quadrant sums into geometric series.  Concerning the
later, we derive the following lemma.

\begin{lemma}  
\label{GeometricSum:lem}
For any $p,q\in\C$ and $K\in\N$, the following is an identity
\be
\label{StdSum}
 \sumQuadK{i}{j} p^i q^j 
     = \fr{1}{2}\fr{p(q^{-1}+2+q)}{(1-p(q+q^{-1})+p^2)} 
        - \fr{1}{2}\lp\fr{1+q}{1-q}\rp 
          \lb \fr{(pq^{-1})^{K+1}}{1-pq^{-1}}
             -\fr{(pq)     ^{K+1}}{1-pq     } \rb. 
\ee
\end{lemma}
\bigskip

\noin {\bf Proof:}\\
The formula follows from iteration of the usual single variable
geometric sum, and algebra. 
\bigskip

We record the following corollary for reference as the expressions appear
many times in the formulas below.

\begin{cor}\label{GeometricSum:cor}
We have the following specializations of lemma~\ref{GeometricSum:lem}
\[
\sumQuadK{m}{n} (\e^{-\l})^m(\e^{-\l\t})^n  
     = 2\e^{-\l}\fOne - 2\e^{-\l(K+1)}f_1^{(K)}(\t,\l)
\]
and
\[
\sumQuadK{n}{m} (\e^{\i\t\l})^n(\e^{\i\l})^m
     = 2\e^{\i\t\l}\fTwo - 2\e^{\i\t\l(K+1)}f_2^{(K)}(\t,\l) 
\]
where the functions $f_1, f_1^{(K)}, f_2, f_2^{(K)}$ are:
\bea
\label{fOne:def}
 \fOne    
  &=& \fr{\cosh^2(\t \l/2)}{1-2\e^{-\l}\cosh(\t \l)+\e^{-2\l}} \\
\label{fOneK:def}
 f_1^{(K)}(\t,\l)
  &=&  \fr{1}{4}\lp\fr{1+\e^{-\l\t}}{1-\e^{-\l\t}}\rp 
       \lb \fr{\e^{ \l\t(K+1)}}{1-\e^{-\l(1-\t)}}
          -\fr{\e^{-\l\t(K+1)}}{1-\e^{-\l(1+\t)}} \rb \\
\label{fTwo:def}
 \fTwo
  &=& \fr{\cos^2(\l/2)}{1-2\e^{\i\t\l}\cos(\l)+\e^{2\i\t\l}} \\
\label{fTwoK:def}
 f_2^{(K)}(\t,\l)
  &=&  \fr{1}{4}\lp\fr{1+\e^{\i\l}}{1-\e^{\i\l}}\rp 
       \lb \fr{\e^{-\i\l(K+1)}}{1-\e^{\i\l(\t-1)}}
          -\fr{\e^{ \i\l(K+1)}}{1-\e^{\i\l(\t+1)}} \rb.
\eea
\end{cor}
\bigskip

We assume $\t$ is in the fundamental region (\ref{FundamentalT}) and
make several observations.  First, we note that all of the
functions given by (\ref{fOne:def})-(\ref{fTwoK:def}) have double
poles at the origin, $\l=0$.  Since $\Im(\t)>0$, neither \fOne\ nor
\fTwo\ have other poles for $\l>0$.  For $\t$ strictly imaginary, the
denominator $1-\exp(-\t\l)$ of $f_1^{(K)}(\t,\l)$ will have isolated
simple zeros.  However, these are balanced by simple zeros of the
difference of bracketed terms in (\ref{fOneK:def}). Thus the product,
in other words $f_1^{(K)}$, has no other singularities for $\l>0$.  A
similar argument shows that $f_2^{(K)}$ is also finite for $\l>0$.
Finally, one may show the bounds
\beas
   |\e^{-\l}\fOne| &<& C_1 \e^{-\l(1-|\Re(\t)|)},\mbox{ and } \\
   |\e^{\i\t\l}\fTwo| &<& C_2 \e^{\i\t\l} 
\eeas
for large $\l$.  As $|\Re(\t)|\le1/2$ and $\Im(\t)>0$, both quantities are 
exponentially decreasing in $\l$.  Similar reasoning shows that
$\e^{-\l(K+1)}f_1^{(K)}(\t,\l)$ and $\e^{\i\t\l(K+1)}f_2^{(K)}(\t,\l)$
are exponentially decreasing in $\l$ {\bf and} $K$.

\section{Eisenstein series}
As mentioned previously, the summation $\tEs$ for $\ReSTwo$ is
absolutely convergent.  We begin by proving our first integral
representation for this case in theorem~\ref{Eisenstein:thm}.  As a
corollary, by restricting $s=n,\ n\ge 3$ we obtain integral
representations of $E_n$.  Further inspection of the integral
representation demonstrates the existence of $\tE_2$.  Elaborating on
theorem~\ref{Eisenstein:thm}, we derive an alternative representation
for $\tEs$ as a contour integral. As a consequence of this second
representation, we prove that \tEs\ admits an analytic continuation in
$s$ as an entire function.  Returning to the analysis of $\tE_2$, we
consider a more general limiting procedure and define \tETwoA\ as the
limit over increasing ``lattice rectangles'' with a fixed aspect ratio
defined by $\al$.  We derive a closed form expression which, when
added to $\tE_2=\tE_2^{(1)}$, gives \tETwoA.  As a corollary, we derive
the relationship between $\tE_2$ and the sum $E_2$ as defined by
Eisenstein.

\subsection{Integral representations}
For the sums $\tEs$ defined by (\ref{HDEn}) we prove
\begin{thm}
\label{Eisenstein:thm}
Given a lattice $\L(\t)$ with ratio $\t$ chosen from the
fundamental region (\ref{FundamentalT}), we have the following
integral representation for $\tEs, \ReSTwo$
\be
\label{EisensteinInt}
\tE_{s}(\t) = 
      \cosPiS \fr{4}{\G(s)}\intRp 
	\l^{s-1}\lp \e^{-\i s\pi/2}\e^{-\l}  f_1(\t,\l)
                  +             \e^{\i\t\l}f_2(\t,\l)\rp \dl
\ee
where \fOne\ and \fTwo\ are given by (\ref{fOne:def}) and
(\ref{fTwo:def}).
\end{thm}
\bigskip

\noin {\bf Proof:}\\
We recognize that, due to the placement of branch cut (\ref{BranchW})
and the symmetric summation conventions (\ref{Grouping}), we have the
following relations between sums over $\LKab{\pm}{\bullet}$ and
$\LKab{\bullet}{\pm}$
\beas
\sumLKab{-}{\bullet} \fr{1}{\w^s} 
  &=& \e^{-\i s\pi} \sumLKab{+}{\bullet} \fr{1}{\w^s} \\
\sumLKab{\bullet}{-} \fr{1}{\w^s} 
  &=& \e^{ \i s\pi} \sumLKab{\bullet}{+} \fr{1}{\w^s}.
\eeas
Therefore we may consider the positive quadrants only and scale the
results by an exponential factor.  Turning to the quadrant
$\LKab{+}{\bullet}$, in place of the isolated singularity of degree
$s$, we substitute the appropriate plane-wave expression from
(\ref{PlaneWave}) to obtain
\beas
\sumLKab{\pm}{\bullet} \fr{1}{(m+n\t)^s}
     &=&(1+\e^{-\i s\pi}) \sumQuadK{m}{n}\fr{1}{(m+n\t)^s} \\
     &=&\fr{(1+\e^{-\i s\pi})}{\G(s)} 
         \intRp \l^{s-1}\sumQuadK{m}{n}\e^{-\l(m+n\t)}\dl \\
     &=&\fr{2(1+\e^{-\i s\pi})}{\G(s)}
         \intRp \l^{s-1}\lp \e^{-\l}\fOne 
             -\e^{-\l(K+1)}f_1^{(K)}(\t,\l)\rp \dl
\eeas
where the last line follows from corollary~\ref{GeometricSum:cor}.
From the statements following this same corollary, we observe that the
two integrands are singular at $\l=0$, and are otherwise finite and
exponentially decreasing in $K$ and $\l\in \R^+$.  In addition, as
$\ReSTwo$, the singularity at the origin is absolutely integrable.
Therefore, one may take the large $K$ limit inside the integral and
compute
\[
 \limK \sumLKab{\pm}{\bullet} \fr{1}{(m+n\t)^{s}}
    =\cosPiS \fr{4}{\G(s)}
     \intRp \l^{s-1}\e^{-\i s\pi/2}\e^{-\l}\fOne \dl.
\]
By a similar analysis we prove that 
\[
\limK \sumLKab{\bullet}{\pm} \fr{1}{(m+n\t)^s} 
  =  \cosPiS \fr{4}{\G(s)}
     \intRp \l^{s-1} \e^{\i\t\l}\fTwo \dl.
\]
Adding these two contributions gives the theorem.  $\dagger$
\bigskip

A few comments are in order.  First, for $s=2j+1$ the term
$\cosPiS=0$.  This is equivalent to the well-known fact that, by
symmetry, the odd Eisenstein series are zero.  Furthermore, by
inspection of (\ref{EisensteinInt}), we see that the absolute
convergence of the Eisenstein series $\tEs,
\ReSTwo$ manifests itself in the behavior of the integrand of
(\ref{EisensteinInt}) near the origin; the factor $\l^{s-1}$ balances
the double poles of $f_1$ and $f_2$ so as to ensure the product is
integrable at $\l=0$.  More careful analysis reveals that the formula
(\ref{EisensteinInt}) is finite even for the conditionally convergent
case $\tE_2$.  The Laurent expansions of the integrands about the
origin are
\beas
 \lim_{\l\to0} \l \e^{-\l}\fr{\cosh^2(\t \l/2)}{1-2\e^{-\l}\cosh(\t \l)+\e^{-2\l}}
    &=& \fr{1}{\l(1-\t^2)} - \fr{(1-2\t^2)\l}{12(1-\t^2)} + \O{\l^3}\\
 \lim_{\l\to0} \l \e^{\i\l\t}\fr{\cos^2(\t \l/2)}{1-2\e^{\i\t\l}\cos(\t \l)+\e^{\i2\t\l}}
    &=& \fr{1}{\l(1-\t^2)} - \fr{(2-\t^2)\l}{12(1-\t^2)} + \O{\l^3}
\eeas
At $s=2$, where the exponential factor $\e^{-\i s\pi/2}=-1$, the
expansions above are to be subtracted.  Hence the integrand of
(\ref{EisensteinInt}) is finite at the origin even in this case.  By a
similar analysis one may show that the $K$-dependent terms also cancel
at the origin.  We have proved:

\begin{cor}\label{EisensteinTwo:cor}
The summation (\ref{HDEn}) converges in the conditionally
convergent case $s=2$, and its value, $\tE_2$, is given by the integral
(\ref{EisensteinInt}).
\end{cor} 

In fact, a great deal more may be said.  The function
\be
\label{EsIntegrand}
   F(s,z) = 
    \e^{-\i s\pi/2}  \e^{-z}\fr{\cosh^2(\t z/2)}{1-2\e^{-z}\cosh(\t z)+\e^{-2z}} 
             + \e^{\i\t z} \fr{\cos^2(z/2)}{1-2\e^{\i\t z}\cos(z)+\e^{2\i\t z}}
\ee
appearing as a factor in the integrand (\ref{EisensteinInt}) has a
singularity at the origin $z=0$, and additional simple poles in the
complex plane at the points
\[
  z \in P = 
   \lsb\left.\pm\fr{2\pi \i}{1\pm \t}m, \pm\fr{2\pi}{1\pm \t}n \right| 
     m,n\in\N\rsb.
\]
We denote the minimum magnitude of all $z\in P$ by $\rho$.  Next,
define the contour $C$ which begins at $\infty+\i y,\ y>0$; runs
parallel to real axis until it intersects the circle centered at the
origin with radius $r$, where $y<r<\rho$; follows this circle
counterclockwise around the origin; and runs back out to $\infty-\i y$,
parallel to the real axis.  We assume that $y>0$ is small enough such
that $C$ encloses only the pole at $z=0$.  Finally, for the function
$z^{s-1}, s\not\in\Z$, situate the branch cut along the positive
real axis such that ($\l\in\R$)
\be 
\label{BranchCut}
\lsb\ba{rcl}
   \Dsp \lim_{y\to 0} (\l+\i y)^{s-1} &=& \l^{s-1} \\
   \Dsp \lim_{y\to 0} (\l-\i y)^{s-1} &=& \e^{2\pi \i (s-1)}\l^{s-1} 
\ea \right.
\ee
With these preliminaries established we prove the following theorem.

\begin{thm}\label{EisensteinCont:thm}
Given a lattice $\L(\t)$ with ratio $\t$ chosen from the
fundamental region (\ref{FundamentalT}), we have the following
contour integral representation for $\tEs,\ \ReSTwo$
\be
\label{EisensteinCont}
\tEs = 2\cosPiS\fr{\G(1-s)\e^{-\i s\pi}}{\i \pi}\int_C z^{s-1}F(s,z)\d z
\ee
where $F(s,z)$ is given by (\ref{EsIntegrand}).
\end{thm}
\bigskip

\noin {\bf Proof:}\\
Consider the contour integral
\[
   \int_C z^{s-1}F(s,z) \d z.
\]
As the integrand is analytic except for the singularity at the origin,
we may deform the contour without altering the value of the integral.
Specifically, shrink the radius of the circle, $r\to 0$, and take
the limit as $y\to 0$ for the two components running parallel to
the real axis.  For $\Re(s)=2+\eps$, we estimate the contribution from
the circular arc
\[
  \lim_{r\to 0} \left| \int_{|z|=r} z^{s-1}F(s,z) \d z \right|
   \le M  \lim_{r\to 0} \int_{|z|=r} r^{\eps-1} |\d z| = 0.
\]
Turning to the components parallel to the real axis, applying the
specification of the branch cut (\ref{BranchCut}) we compute
\beas
  \lim_{y\to0} \int_{\infty+\i y}^{\eps_r+\i y} z^{s-1}F(s,z)\d z 
    &=& -\intRp \l^{s-1}F(s,\l) \dl \\
  \lim_{y\to0} \int_{\eps_r-\i y}^{\infty-\i y} z^{s-1}F(s,z)\d z 
    &=& \e^{2\pi \i(s-1)}\intRp \l^{s-1}F(s,\l) \dl.
\eeas
We recognize the integrals of theorem~\ref{Eisenstein:thm},
make use of the identity
\[
  \G(s)\G(1-s) = \fr{\pi}{\sin(\pi s)},
\]
and obtain
\beas
    \int_C z^{s-1}F(s,z) \d z 
   &=& (\e^{2\pi \i(s-1)}-1)\intRp \l^{s-1}F(s,\l) \dl \\
   &=& \fr{2\i \e^{\pi \i s}\sin(\pi s) \G(s)}{4\cosPiS} \tEs \\
   &=& \fr{\i\pi}{2\cosPiS \e^{-\i s\pi}\G(1-s)} \tEs.
\eeas
The result (\ref{EisensteinCont}) follows from algebra. $\dagger$
\bigskip

Several corollaries follow from theorem~\ref{EisensteinCont:thm}.  We
note here only the most immediate

\begin{cor}\label{EsEntire:cor}
The sums $\tEs$, defined by (\ref{HDEn}) for \ReSTwo, admit an
analytic continuation to $s\in\C$ as an entire function.  This
continuation is given by the contour integral representation
(\ref{EisensteinCont}).
\end{cor}
\bigskip

\noin {\bf Proof:}\\ 
The contour integral appearing in (\ref{EisensteinCont}) defines an
analytic function of $s$ which is never singular.  The same may be
said for the cosine factor.  Thus the only possible singularities
would arise from the factor $\G(1-s)$ which has simple poles for
$s\in\N$.  From the definition (\ref{HDEn}) we know that \tEs\ is
finite for \ReSTwo\ (the apparent singularities in
(\ref{EisensteinCont}) in this case are balanced by zeros of the
cosine term, the contour integral, or both.)  Therefore, we need only
verify the finite existence of $\tE_1$ and $\tE_2$.  We have argued
above that $\tE_2$ is finite.  Finally, for $s=n=1$, the pole in the
Gamma function is balanced by the simple zero of the cosine
factor. $\dagger$
\bigskip

By inspection of (\ref{EisensteinCont}), we find that $\tEs$ has
simple zeros for $s=1-2j,\ j\ge1$.  Thus the continuation respects the
symmetric limiting process (\ref{HDEn}).  We anticipate further
results concerning the evaluation of the contour integral
(\ref{EisensteinCont}) via residue methods. However, we have not
carried out this analysis at the time of writing this paper.  Finally,
as mentioned previously, the exact form of the integrands
(\ref{EisensteinInt}) and (\ref{EisensteinCont}) reflect our summation
convention with respect to grouping of summands and placement of the
branch-cut.  Regarding the latter, similar formulas arise if the
branch-cut is situated along any of the lattice diagonals; the effect
is to redistribute factors of $\exp(\i\pi s/2)$ between the two
functions $f_1$ and $f_2$.  There are, perhaps, additional treatments
of the branch-cut that could yield relatively simple expressions.
However, a simple integral expression valid for placement of the cut
along an arbitrary ray in the complex plane appears to be intractable.

\subsection{Aspect ratio correction}
Although we have demonstrated that $\tE_2$ is finite, until this point
its relation to Eisenstein's definition of this series is unclear.  We
make this explicit in the present section.

We wish to formalize summation over lattice rectangles with a fixed
aspect ratio.  Given $\al\in(0,\infty)$ we define the sum
\be
\label{HDEab} 
\tETwoA = 
    \limK \sum_{|m|\le\alK}\sum_{|n|\le K} \fr{1}{(m+n\t)^2}, 
\ee
where $\lfloor x\rfloor$ denotes greatest integer less than or equal
to $x$. Again, although expected, the existence of the limit
(\ref{HDEab}) is not {\em a priori} guaranteed but will follow in the
course of our analysis.

Clearly, $\tE_2^{(1)}=\tE_2$ defined in (\ref{HDEn}).  More generally,
we write
\be
\label{SumSplitting}
\tETwoA = \tE_2 + \Delta(\al,\t)
\ee
where the value of $\tE_2$ may be computed via the integral expression
(\ref{EisensteinInt}).  Analysis similar to that used to prove
(\ref{EisensteinInt}) may be employed to compute a closed form
expression for $\Delta(\al,\t)$.

\begin{thm}
\label{AspectRatio:thm}
For a fixed aspect ratio determined by $\al\in(0,\infty)$, the limit
$\tETwoA$ specified in (\ref{HDEab}) exists.  Furthermore, when
written in the form (\ref{SumSplitting}), the aspect ratio dependence
is given by
\be 
\label{AspectRatio}
\Delta(\al,\t) = -\fr{4\i}{\t}
               \lp\arctan\lp \i\t\rp - \arctan\lp\fr{\i\t}{\al}\rp\rp
\ee
\end{thm}
\bigskip

\noin {\bf Proof:}\\
Assume $\al\ge1$.  We write the limit (\ref{HDEab})
\beas
 \tETwoA 
   &=& \tE_2 + \Delta(\al,\t),\\
\Delta(\al,\t)  
   &=& 2\limK
     \sum_{m=K+1}^{\alK}\sum_{n=-K}^{n=K} \fr{1}{(m+n\t)^2}.
\eeas
The contribution from the sum over lattice points $-\lfloor\al\cdot
K\rfloor \le m \le -K-1$ is accounted for by the factor of two
multiplying the sum in the final line.  We represent the poles using
the $\LKab{+}{\bullet}$ plane-wave expansion (\ref{PlaneWave}).
\bea
 \Delta(\al,\t) 
   &=& 2\limK
     \sum_{m=K+1}^{\alK} 
     \sum_{n=-K}^{n=K} \intRp \l \e^{-\l(m+n\t)} \dl \nn\\
   &=& 2\limK
     \intRp \l 
     \lp\fr{\e^{\l\t K}-\e^{-\l\t(K+1)}}{1-\e^{-\l\t}}\rp
     \lp\fr{\e^{-\l(K+1)}-\e^{-\l\alK}}{1-\e^{-\l}}\rp  \dl \nn\\
   &=& 2\limK
     \intRp \l 
     \lp\fr{\e^{\l\t}-\e^{-\l\t(1+1/K)}}{1-\e^{-\l\t/K}}\rp
     \lp\fr{\e^{-\l(1+1/K)}-\e^{-\l\alK/K}}{1-\e^{-\l/K}}\rp  
     \fr{1}{K^2}\dl \nn \\
\label{AlphaInt}
   &=& \fr{2}{\t}
     \intRp \fr{(\e^{\l\t}-\e^{-\l\t})(\e^{-\l}-\e^{-\l\al})}{\l} \dl.
\eea
Where the argument which justifies taking the large $K$ limit inside
the integral runs along the same lines as in the proof of
theorem~\ref{Eisenstein:thm}.  This last integral (\ref{AlphaInt}) may
be evaluated in closed form using the formula
\be
\label{IntegralFormula}
  \intRp \e^{-\beta x}\sin(\delta x) \fr{1}{x}\d x
     = \arctan\lp\fr{\delta}{\beta}\rp, 
\ee
which holds for $\Re(\beta)>|\Im(\delta)|$ (see \cite{gr94}, 3.944.5).
Taking care to write (\ref{AlphaInt}) as the difference of two
integrals of the form (\ref{IntegralFormula}), and performing algebra
gives the expression for $\Delta(\al,\t),\ \al\ge 1$ in
(\ref{AspectRatio}).  

For $\al<1$ the lattice rectangle is such that the longer side is in
the $\t$-direction.  As written, equations (\ref{HDEab}) and
(\ref{SumSplitting}) suggest that this rectangle is inscribed in a
lattice square of size $K$, and to compute $\Delta(\al,\t)$, one should
subtract the extra contributions exterior to the rectangle but
interior to the square.  The problem with this approach is that, for
arbitrary $\t$ and $\al$, one would have to keep track of the
quadrants in which these points lie.  In lieu of this, for $\al<1$ we
rescale the limits in (\ref{HDEab})
\[
\tETwoA = \limK \sum_{|m|\le K}\sum_{|n|\le \betaK} \fr{1}{(m+n\t)^2}.
\]
Informally, this has the ``effect'' of inscribing the square in the
rectangle and motivates computing the contributions from the points in
the difference using the plane-wave formulas appropriate for
$\LKab{\bullet}{\pm}$.  Arguing as above and using the integral
identity (\ref{IntegralFormula}) we compute
\bea
\Delta(\al,\t) 
 &=&-2\limK
     \sum_{n=K+1}^{\betaK} 
     \sum_{m=-K}^{m=K} \intRp \l \e^{\i\l(m+n\t)} \dl \nn \\
 &=& \fr{2}{\t}
     \intRp \fr{(\e^{-\i\l}-\e^{\i\l})(\e^{\i\l\t}-\e^{\i\l\t/\al})}{\l} \dl \nn\\
\label{BetaInt}
 &=& -\fr{4\i}{\t}\lp\arctan\lp-\fr{1}{\i\t}\rp-\arctan\lp-\fr{\al}{\i\t}\rp\rp.
\eea

Although perhaps not obvious at first glance, the formula
(\ref{BetaInt}) is the same as the formula for $\Delta(\al,\t)$ in
(\ref{AspectRatio}).  Using standard trig identities we have
\be
\label{ArcTanIdent}
\lim_{z'\to z} \lp\arctan(z) - \arctan\lp-\fr{1}{z'}\rp\rp
    = \lim_{z'\to z} \arctan\lp\fr{zz'+1}{z'-z}\rp
    = \pm \fr{\pi}{2}.
\ee
Furthermore, for $z=\i\t$ or $z=\i\t/\al$ with $\t$ satisfying
(\ref{FundamentalT}) and $\al>0$, we find that we should choose the
minus sign in (\ref{ArcTanIdent}).  Using this identity and algebra we
observe that (\ref{BetaInt}) is equal to (\ref{AspectRatio}), thus
(\ref{AspectRatio}) holds for all $\al>0$.  $\dagger$
\bigskip

We use this theorem to find the connection between our summation and
Eisenstein's.  Starting from Eisenstein's summation convention we
obtain
\bea
 E_2 &=& \sumE \frac{1}{(m + n\t)^{2}} \nn \\
     &=&\lim_{N\to\infty}\sum_{n=-N}^N 
        \lim_{\al\to\infty}
        \sum_{m=-\lfloor \al\cdot N\rfloor}^{m=-\lfloor \al\cdot N\rfloor}
        \fr{1}{(m+n\t)^2} \nn \\
     &=&\lim_{\al\to\infty}
        \lim_{N\to\infty}\sum_{n=-N}^N 
        \sum_{m=-\lfloor \al\cdot N\rfloor}^{m= \lfloor \al\cdot N\rfloor}
        \fr{1}{(m+n\t)^2} \nn \\
\label{EisensteinToHD}
     &=&\tE_2 -\fr{4\i}{\t}\arctan(\i\t).
\eea
Standard estimates justify commuting the $N$ and $\al$ limits between
the second and third lines, and we used (\ref{AspectRatio}) to compute
this limit.  As $\t\not=0$, (\ref{EisensteinToHD}) shows that, for
finite $\t$, the value of our sum is always different from
Eisenstein's.  In the limit $\t\to \i\infty,\ |\Re(\t)|\le 1/2$, both
summations are equal and presumably converge to $2\z_R(2)=\pi^2/3$.

In a similar vein, we compute the difference between taking
Eisenstein's limit and ``its reverse''.  Arguing as above we have that
\bea
\sumNotE \fr{1}{(m+n\t)^2}
  &=& \tE_2 +\lim_{\al \to 0}\Delta(\al,\t) \nn\\
\label{NotEisensteinToHD}
  &=& \tE_2 -\fr{4\i}{\t}\lp\arctan\lp\fr{\i}{\t}\rp + \fr{\pi}{2}\rp
\eea
Taking the difference between (\ref{EisensteinToHD}) and
(\ref{NotEisensteinToHD}) and we obtain
\[
\lp \sumE - \sumNotE \rp \fr{1}{(m+n\t)^2}
  = \fr{2\pi \i}{\t}.
\]
For a different proof of this fact see Walker, \cite{walker96}.

Finally, in the case of the square lattice, $\t=\i$, we observe that
$f_1(\i,\l)=f_2(\i,\l)$.  Collecting factors in (\ref{EisensteinInt})
and performing algebra we find that
\[
  \tEs(\i) = \e^{-\i s\pi/4}
            \cosPiS \cos\lp\fr{\pi}{4}s\rp \fr{8}{\G(s)}
            \intRp \l^{s-1}\e^{-\l}f_1(\i,\l) \dl.
\]
Corresponding to the added $\pi/2$-symmetry of the square lattice, the
product of cosines causes the sum to vanish for $s=n\ge 3$, $n$ not a
multiple of four.  Similarly, as the representation
(\ref{EisensteinInt}) is valid for $n=2$, we conclude $\tE_2(\i)=0$.
Substituting this value into (\ref{SumSplitting}), and taking the
large $\al$ (small $\al$) limits we obtain
\beas
  \sumE    \fr{1}{(m+\i n)^2} &=& \pi \\
  \sumNotE \fr{1}{(m+\i n)^2} &=&-\pi,
\eeas
well-known identities in the fast multipole community (see, for
example, \cite{greengard88}).

\section{Weierstrass elliptic functions}
Our derivations of integral formulas for Weierstrass' elliptic
functions proceed in much the same manner as above.  As with the
Eisenstein series, we begin with the definitions of the functions as
sums over the lattice (\ref{WeierstrassP}) and (\ref{WeierstrassZ}).
Next, we group terms in the sum as in (\ref{Grouping}), substitute the
appropriate plane-wave expansions from lemma~\ref{PlaneWave:lem}, and
sum the resulting geometric series using the identities of
corollary~\ref{GeometricSum:cor}.

As a preliminary note, the integral representations for $\P(z,\t)$ and
$\z(z,\t)$ which we derive in theorem~\ref{Weierstrass:thm} are not
valid for all $z\in\C$, but rather have a finite domain of validity.
This is a consequence of the way in which we group terms. More
precisely, we require that $z\in D(\t)$ defined by
\be
\label{FundamentalZ}
D(\t) = \{z |\ \Re(-1 \pm z \pm \t)<0,\ \Im(\t \pm z )>0\}.
\ee
As $\t$ is in the region (\ref{FundamentalT}), one may verify that
$D(\t)$ is an open set containing the origin.  However, we point out
that $D(\t)$ may not contain the fundamental period parallelogram of
the lattice.  For example, the standard hexagonal lattice has
generators $(1,\t)=(1,1/2+\i\sqrt{3}/2)$.  Thus a corner of $\L_0$ is
given by the point $z_0=1/2+\t/2=3/4+\i\sqrt{3}/4$.  However,
$\Re(-1+z_0+\tau)=1/4>0$, violating the first inequality in
(\ref{FundamentalZ}).

With this aside, we prove the following.
\begin{thm}\label{Weierstrass:thm}
Assume $\L=\L(\mu,\nu)$ is an arbitrary complex lattice with
generators chosen such that $\t=\nu/\mu$ is in the fundamental region
(\ref{FundamentalT}), and that the complex number $z$ is in the domain
$D(\t)$ defined by the inequalities (\ref{FundamentalZ}).  We have the
following integral expressions for the inhomogeneous elliptic functions
$\P(z,\t)$, and $\z(z,\t)$:
\be
\label{WeierstrassPInt}
 \P(z,\t) = \Dsp \fr{1}{z^2} + 8\intRp 
	\l\lb \e^{-\l}    \sinh^2\lp\fr{z\l}{2}\rp f_1(\l,\t)
             +\e^{\i\t \l}\sin^2 \lp\fr{z\l}{2}\rp f_2(\l,\t)\rb \dl,
\ee
and,
\be
\label{WeierstrassZInt}
 \z(z,\t) = \Dsp \fr{1}{z} + 4\intRp 
	  \lb \e^{-\l}    (z\l-\sinh(z\l))f_1(\l,\t)
             -\e^{\i\t \l}(z\l-\sin(z\l) )f_2(\l,\t)\rb \dl,
\ee
where the functions $f_1, f_2$ are defined by (\ref{fOne:def}) and
(\ref{fTwo:def}).  We may evaluate the homogeneous functions,
$\P(x|\mu,\nu)$ and $\z(x|\mu,\nu)$, via the appropriate scaling
relations (\ref{Rescalings}) and
(\ref{WeierstrassPInt})-(\ref{WeierstrassZInt}) with the proviso that
$x/\mu=z\in D(\t)$.
\end{thm}
\bigskip

\noin {\bf Proof:}\\
In computing the integral representation for the \P\ function we will
group terms of the sum (\ref{WeierstrassP}) as in (\ref{Grouping}).
As in the computation of the Eisenstein sums, we wish to combine the
contributions from the sums over $\LKab{\pm}{\bullet}$.  We compute
\beas
\sumLKab{-}{\bullet}  \lp \fr{1}{(z-\w)^2}-\fr{1}{\w^2} \rp
  &=& \sumQuadK{m}{n} \lp \fr{1}{(z+m+n\t)^2}-\fr{1}{(m+n\t)^2}\rp \\
\sumLKab{+}{\bullet}  \lp \fr{1}{(z-\w)^2}-\fr{1}{\w^2} \rp
  &=& \sumQuadK{m}{n} \lp \fr{1}{(-z+m+n\t)^2}-\fr{1}{(m+n\t)^2}\rp.
\eeas
Therefore the contributions from both quadrants may be expressed as a
single sum over $\LKab{+}{\bullet}$ of a modified summand.
Furthermore, under the assumption $z\in D$, all of the poles in this
sum may be expressed using the $\Re(\w)>0$ plane-wave expansion from
(\ref{PlaneWave}).
\bea
\sumLKab{\pm}{\bullet} \lp \fr{1}{(z-\w)^2}-\fr{1}{\w^2} \rp
 &=& \sumQuadK{m}{n} 
     \lp\fr{1}{(-z+m+n\t)^2}-\fr{2}{(m+n\t)^2}+\fr{1}{(z+m+n\t)^2}\rp \nn\\
 &=& \intRp \l(\e^{\l z}-2+\e^{-\l z})\lp 
            \sumQuadK{m}{n} \e^{-\l(m+n\t)}\rp \dl \nn \\
\label{WeierstrassPIntKm}
    &=& 8 \intRp \l\sinh^2\lp\fr{z\l}{2}\rp \lp \e^{-\l}\fOne -  
               \e^{-\l(K+1)}f_1^{(K)}(\t,\l) \rp \dl.
\eea
Arguing as before, we find the large $K$ limit of the $K$-dependent
term to be zero.  We compute the contribution from the terms in the
quadrants $\LKab{\bullet}{\pm}$ in an analogous manner.  Adding this
result to (\ref{WeierstrassPIntKm}) gives (\ref{WeierstrassPInt}).

The derivation of the expression for the \z\ function is similar.  In
brief, the sum over the quadrants $\LKab{\pm}{\bullet}$ may again be
expressed as a sum over the single quadrant $\LKab{+}{\bullet}$ in
which we substitute the appropriate plane-wave expansion.  Thus,
\bea
\sum_{\LKab{\pm}{\bullet}} \lp \fr{1}{(z-\w)}+\fr{1}{\w}+\fr{z}{\w^2} \rp
  &=& \sumQuadK{m}{n} 
      \lp -\fr{1}{-z+m+n\t}+\fr{2z}{(m+n\t)^2}+\fr{1}{z+m+n\t}\rp  \nn\\
  &=& \intRp (\e^{\l z}+2z\l-\e^{-\l z})\lp 
             \sumQuadK{m}{n} \e^{-\l(m+n\t)}\rp \dl \nn \\
\label{WeierstrassZIntKm}
  &=& 4 \intRp (z\l - \sinh(z\l)) \lp \e^{-\l}f_1(\l,\t)
      -\e^{-\l(K+1)}f_1^{(K)}(\l,\t) \rp \dl.
\eea
As before the $K$-dependent term goes to zero in the limit.  The
analogous sums over $\LKab{\bullet}{\pm}$ give the other half of the
expression (\ref{WeierstrassZInt}). $\dagger$
\bigskip

\noin Remarks: 
\begin{enumerate}
\item As an alternative to the above derivation of
theorem~\ref{Weierstrass:thm}, we recall that the Eisenstein series
appear as coefficients in the Laurent expansion for Weierstrass' \P\
function
\beas
 \P(z)  &=&\fr{1}{z^2} + \sumLat \lp \fr{1}{(z-\w)^2}-\fr{1}{\w^2}\rp \\
        &=&\fr{1}{z^2} + \sumLat \sumTay (n-1)z^{n-2}\fr{1}{\w^n}\\
        &=&\fr{1}{z^2} + \sumTay (n-1) E_n z^{n-2}.
\eeas
Substituting the integral representations (\ref{EisensteinInt}) for
the $E_n$, the Taylor series may be summed explicitly inside the
integrand.  The formula (\ref{WeierstrassPInt}) above follows after
algebraic simplification.  Furthermore, the expression
(\ref{WeierstrassZInt}) for the \z\ function follows from
anti-differentiation of (\ref{WeierstrassPInt}).

\item As with the series $\tE_2$, and its dependence on aspect ratios 
derived in theorem~\ref{AspectRatio:thm}, the slowly decaying terms of
the sums defining \P\ and \z\ manifest themselves at the origin in the
integral representations (equations (\ref{WeierstrassP}),
(\ref{WeierstrassZ}) and (\ref{WeierstrassPInt}),
(\ref{WeierstrassZInt}) respectively).  In the integral
representations, we observe that Weierstrass' ``correction'' terms are
arranged in such a way as to create third order zeros at $\l=0$, which
appropriately balance the second order poles from $f_1$ and $f_2$.
\end{enumerate}

\section{Conclusion}
We conclude with a brief discussion of our results in relation to
previous research in this field.  To the best of our knowledge, there
is no analog to the integral expressions for the \P\ and \z\ functions
(\ref{WeierstrassPInt}) and (\ref{WeierstrassZInt}).  The possibility
of developing numerical routines for evaluation of these functions
based on these representations deserves further study.  We observe
that the integrands are not extremely oscillatory, and decay
exponentially.  Thus $N$-point Gauss-Laguerre quadrature rules will
converge rapidly in $N$. As one drawback, there is the perhaps awkward
domain of validity in $z$.  However, it may be that symmetries
of the \P\ and \z\ functions imply that it is sufficient to
evaluate them over domains that are smaller than the
fundamental period parallelogram.  Furthermore, at least for the \P\
function, there exists the following closed-form Fourier-like
expansion \cite{walker78}
\be
\label{WeierstrassPFourier}
  \P(z,\t) = 
   -2\lp \fr{1}{6} + \sum_{n=1}^\infty \fr{1}{\sin^2(n\pi\t)}\rp
   + \fr{\pi^2}{\sin^2(\pi z)}
   -8\pi^2 \lp \sum_{n=1}^{\infty} 
     \fr{n \e^{2\pi \i\t n}}{1-\e^{2\pi \i\t n}} \cos(2\pi n z)\rp .
\ee
Both summands in (\ref{WeierstrassPFourier}) are exponentially
decreasing and the sums converge rapidly---stiff competition from a
numerical perspective.  Nevertheless, we have not fully explored the
relative merits of this approach over the plane-wave representation
(\ref{WeierstrassPInt}).  In addition, the integral representations
may have further analytic implications.

Turning to the representations for the Eisenstein series, the
existence of $\tE_2$ (corollary~\ref{EisensteinTwo:cor}) is not
unexpected.  In addition to the original finiteness proofs given by
Eisenstein, many years prior to this present work, Walker derived
the remarkable formula for the conditionally convergent series (see
\cite{walker78})
\[
  \limK \sum_{0<m^2+n^2\le K^2} \fr{1}{(m+n\t)^2} 
    = \fr{-2\pi}{1-\i\t} - 4\pi \i \fr{\eta'(\t)}{\eta(\t)}, 
\]
where the Dedekind $\eta(\t)$-function with $\Im(\t)>0$ is defined by
\[
 \eta(\t) = \e^{\pi \i \t/12}\prod_{n=1}^\infty (1-\e^{2\pi \i\t n}).
\]
We also note that a different treatment, initiated by Hecke,
has become a standard approach to resolving convergence and
transformation properties of $E_2$ \cite{schoeneberg74}.

As indicated by theorem~\ref{EisensteinCont:thm}, our expressions are
quite general, and have broad implications.  We note that Riemann
demonstrated both the functional equation satisfied by $\z_R(s)$, and
the evaluation of $\z_R(-2n+1)$ (and, via the functional equation,
$\z_R(2n)$) in terms of Bernoulli numbers using the ``version'' of
theorem~\ref{EisensteinCont:thm} appropriate for his zeta function.
Similarly, we anticipate that a residue argument will give the
evaluation of $\tE_n=E_n$ in terms of multiple Bernoulli numbers.  For
an alternative treatment of Eisenstein series for negative even
integers using Hecke convergence factors see the recent work of
Pribitkin \cite{pribitkin00}.  The functional equation satisfied by
the continuation of $\tEs$ is more elusive.  We are currently pursuing
this and hope to report our results in the future.

Finally, there is a possibility that representations of the form
(\ref{EisensteinInt}) may exist for certain Dirichlet series
\beas
  G(s,\chi) &=& \sumLat \fr{\chi(\w)}{|\w|^{2s}} \\
  \chi(\w)  &=& \exp(\i(m\mu\al + n\nu\beta)),
\eeas
for $\al, \beta \in \R$. (These are called ``Kronecker series'' in
\cite{weil76}, Chapter VIII.)  A detailed discussion of the
convergence of these series is given in \cite{bbp98}.  We note that
Laplace-Mellin techniques have been employed frequently in the
analysis of such series.  The approach up until now has been to think
of
\[
  |\w_{m,n}|^2 = |m\mu + n\nu|^2 = Q(m,n)
\]
as defining a positive definite quadratic form taking $m$ and $n$ as
arguments.  Treating this form as ``indivisible'', one may use the
$\Re(\w)>0$ plane-wave formula in lemma~\ref{PlaneWave:lem} and obtain
the integral representation
\[
\sumLat \fr{\chi(\w)}{|\w|^{2s}} 
   = \fr{1}{\G(s)} \intRp \l^{s-1}  \sumLat \chi(\w)\e^{-\l Q(m,n)}\dl.
\]
The analysis then proceeds via $\theta$ functions.

Our approach would be different.  ``Plane-wave-like'' representations
exist for the function $f(x,y)=\sqrt{x^2+y^2}$.  Formally, one may
take the true plane-wave expressions for $f(x,y,z)=\sqrt{x^2+y^2+z^2}$
derived in \cite{cgr99}, and set $z=0$.  The result is a 2D
integral --- as opposed to the Eisenstein case analyzed above where one
complex dimension (two real) collapses into a 1D integral.  However,
the critical element of this representation is that the exponential
function in the ``plane-wave'' representation is linear in $m$ and
$n$.  Again, as a consequence, we observe that the summation under the
integrand becomes a 2D geometric series.  We are considering this as a
possible direction for future research.
\bigskip

\noin{\bf Acknowledgments.}  The authors thank Peter Walker at the
University of Sharjah, United Arab Emirates, and Brad Alpert at
NIST/Boulder for frequent helpful discussions during the writing of
this paper.

\newpage

\newpage
\bibliography{/home/andrewd/Bibliography/master}
\bibliographystyle{abbrv}

\end{document}